\documentstyle[amscd]{amsart}
\makeatletter
\renewcommand{\subsection}{\@startsection{subsection}{2}{\z@}%
{\baselineskip}{0.5\baselineskip}{\defaultfont\bf}}
\makeatother
\newsymbol\varnothing 203F
\newtheorem{thm}{Theorem}[section]
\newtheorem{pro}[thm]{Proposition}
\newtheorem{lem}[thm]{Lemma}
\theoremstyle{definition}
\newtheorem{defn}{Definition}[section]
\numberwithin{equation}{section}
\def\dj{d\kern-.30em\raise1.25ex\vbox{\hrule width .3em height .03em}}
\def\Dj{D\kern-.75em\raise0.75ex\vbox{\hrule width .3em height .03em}
\kern.03em}
\newcommand{\hor}{\mbox{\family{euf}\shape{n}\selectfont hor}}
\def\k{\kappa}
\def\bla#1{$(${\it #1\/{}}$)$}
\def\sstar{{\raise0.2ex\hbox{$\scriptstyle\star$}}}
\def\btwn{\gamma}
\def\im{\mathrm{im}}
\def\tr{\mathrm{tr}}
\def\bimh#1{\cal{F}_{#1}}
\def\map#1{\varrho_{#1}}
\def\pre#1{#1_{\star}}
\def\amap#1{\sharp_{#1}}
\def\bim#1{\cal{E}_{#1}}
\def\Mor{\mathrm{Mor}}
\def\MoR{\mathrm{M}}
\def\M{\bar{\mathrm{M}}}
\def\H{\bar{\mathrm{H}}}
\def\Hom{\mathrm{Hom}}
\def\HoM{\mathrm{H}}
\def\R#1{R(#1)}

\def\1{\varnothing}
\def\I{I}

\newcommand{\id}{\mathrm{id}}
\newcommand{\e}{\epsilon}
\newcommand{\Sum}{\displaystyle{\sum}}
\begin{document}
\title[Quantum Principal Bundles]{Quantum Principal Bundles $\&$
Tannaka-Krein Duality Theory}
\author{Mi\'co {\Dj}ur{\Dj}evi\'c}
\address{Instituto de Matematicas, UNAM, Area de la Investigacion
Cientifica, Circuito Exterior, Ciudad Universitaria, M\'exico DF, CP
04510, MEXICO}
\maketitle
\begin{abstract}
The structure of quantum principal bundles is studied, from the viewpoint
of Tannaka-Krein duality theory. It is shown that if the structure
quantum group is compact, principal $G$-bundles over a quantum space $M$
are in a natural correspondence with certain contravariant functors
defined on the category of finite-dimensional unitary representations of
$G$, with the values in the category of finite projective bimodules over
a *-algebra representing the base space.
\end{abstract}
\renewcommand{\thepage}{}
\section{Introduction}
The aim of this paper is to study the internal structure of quantum
principal bundles possessing compact structure groups, from the point of
view of Tannaka-Krein duality theory. This analysis can be understood as
a non-commutative generalization of the classical incorporation of the
structure of principal bundles in the language of corresponding
associated vector bundles.

The structure of a compact matrix quantum group $G$ is completely
encoded \cite{W1,W2} in the category $\R{G}$ of finite-dimensional
unitary representations of
$G$ (corresponding to $G$, in the framework of Tannaka-Krein duality).
The main result of this paper is that there exists a natural
correspondence between quantum principal $G$-bundles $P$ over a quantum
space $M$, and appropriate contravariant functors $\tau_P$ on $\R{G}$,
with values in the category $\Phi_M$ of finite projective bimodules
over the *-algebra $\cal{V}$ representing the base $M$.
Explicitly, these bimodules consist of intertwiners between
representations from $\R{G}$ and the right action
of $G$ on $P$. They are interpretable as quantum
counterparts of associated vector bundles.

In non-commutative differential geometry \cite{C}, vector bundles can be
naturally viewed as one-sided finite-projective modules. Such a concept
of a vector bundle is especially natural in considerations involving
cyclic cohomology. In particular, the standard algebraic $K$-theory
naturally enters the game. However, associated vector bundles
introduced in this paper are much more rigid objects, in particular they
always appear as {\it bimodules}.

The paper is organized as follows. The next section describes a
construction of the functor $\tau_P$, starting from a quantum principal
bundle $P$. In Section~3 the inverse construction is presented,
which reconstructs the bundle starting from a contravariant
functor $\tau\colon\R{G}\rightarrow\Phi_M$
possessing necessary additional properties. In Section~4 some
examples are considered.

It is worth noticing that even if $M$ and $G$ are classical the bundle
$P$ may be a quantum object. In particular, an interesting purely
quantum phenomena is the nontriviality of the classification for
bundles over a $1$-point set.

The important point is that functors $\tau$ do not see the ``concrete''
aspect of the category $\R{G}$, given by specifying carrier unitary
spaces of representations. This means that the classification problem
for groups having ``abstractly'' the same categories of representations
will be the same.

\section{Structural Analysis of Quantum Principal Bundles}

Let $G$ be a compact matrix quantum group \cite{W1}, represented by a
Hopf *-algebra $\cal{A}$, consisting of ``polynomial functions'' on $G$.
Let $\phi\colon\cal{A}\rightarrow\cal{A}\otimes\cal{A}$,
$\e\colon\cal{A}\rightarrow\Bbb{C}$ and
$\k\colon\cal{A}\rightarrow\cal{A}$ be the coproduct, counit and the
antipode map.

In the framework of Tannaka-Krein duality theory \cite{W2}, the group
$G$ is represented by a complete concrete monoidal $W$-category
$\R{G}$. The objects of $\R{G}$ are unitary representations  of $G$
in finite-dimensional unitary spaces. Morphisms in $\R{G}$
are intertwiners between representations.

For each $u\in\R{G}$ we shall denote by $H_u$ the corresponding
carrier space (so that $u$ is understandable as a map $u\colon H_u
\rightarrow H_u\otimes\cal{A}$). We shall denote by $\oplus$, $\times$
the sum and the product in $\R{G}$. In particular,
$$ H_{u\oplus v}=H_u\oplus H_v\qquad H_{u\times v}=H_u\otimes H_v.$$
We shall denote by $\1$ the trivial representation of $G$, acting in
$H_{\!\1}=\Bbb{C}$. In what follows we may assume that spaces $H_u$ belong
to a given model set of unitary spaces, so that $\R{G}$ is small.
\renewcommand{\thepage}{\arabic{page}}

Let us consider a quantum space $M$, represented by a unital
*-algebra $\cal{V}$. Let $P=(\cal{B},i,F)$ be a quantum principal
$G$-bundle \cite{D2} over $M$. Here $\cal{B}$ is a *-algebra consisting
of appropriate ``functions'' on $P$, while
$F\colon\cal{B}\rightarrow\cal{B}\otimes\cal{A}$ and
$i\colon\cal{V}\rightarrow\cal{B}$ are
unital *-homomorphisms playing the role of the ``dualized''
right action of $G$ on $P$, and the dualized projection of $P$ on $M$
respectively. We have
$$ (F\otimes\id)F=(\id\otimes\phi)F\qquad (\id\otimes\e)F=\id\qquad
i(\cal{V})=\Bigl\{\,b\in\cal{B}; F(b)=b\otimes 1\,\Bigr\}.$$

Let us consider the
intertwiner $\cal{V}$-bimodules
$$\bim{u}=\Mor(u,F).$$
The elements of these bimodules are interpretable as smooth sections
of corresponding associated vector bundles.

{}From the analysis of \cite{D2}--Appendix A it follows that $\cal{E}_u$
are (both side) finite and projective. Each morphism
$\lambda\in\Mor(u,v)$ induces a bimodule homomorphism
$\map{uv}(\lambda)\colon\bim{v}\rightarrow\bim{u}$, via the composition
of intertwiners. In other words, we have a natural
system of linear maps
$\map{uv}\colon\Mor(u,v)\rightarrow\Hom(\bim{v},\bim{u})$. In what
follows we shall also use a simpler alternative notation and write
$$\pre{\lambda}=\map{uv}(\lambda),$$
when the domain/codomain is clear from the context.

We have
\begin{equation}\label{contr}
\pre{\lambda}\pre{\rho}=\pre{[\rho\lambda]}
\end{equation}
for $\lambda\in\Mor(u,v)$ and $\rho\in\Mor(v,w)$.

For each $u\in\R{G}$ let $\bar{u}=u^c$ be the corresponding conjugate
representation, identified with the contragradient representation.
By definition,
$$ H_{\bar{u}}=H_u^*\qquad u^c=(\id\otimes\k)(u)^\top.$$

Let $C_u\colon H_u\rightarrow H_u$ be the canonical intertwiner between $u$ and
$u^{c\!c}$. This is a strictly positive map, normalized such that
$\tr(C_u)=\tr(C_u^{-1})$.
Explicitly, $C_u$ can be described as follows. Let $(f_z)_{z\in\Bbb{C}}$
be a canonical holomorphic family of linear multiplicative functionals
$f_z\colon\cal{A}\rightarrow\Bbb{C}$, describing modular properties of
the Haar measure \cite{W1}. Then we have
$$ C_u=(\id\otimes f_1)u.$$
In particular, it follows that
$$ C_u\otimes C_v=C_{u\times v}\qquad \varphi C_u=C_v\varphi $$
for each $\varphi\in\Mor(u,v)$ and $u,v\in\R{G}$. The
natural scalar product in $H_u^*$ is given by
$$(f,g)=(j_u^{-1}(g),C_uj_u^{-1}(f))$$
where $j_u\colon H_u\rightarrow H_u^*$ is the canonical antilinear map
(induced by the scalar product in $H_u$). The representation
$\bar{u}$ is unitary, with respect to the above scalar product.
Furthermore, we have
$$C_u=j_{\bar{u}}j_u$$
for each $u\in\R{G}$.

In what follows, a natural identification
$(W\otimes V)^*=V^*\otimes W^*$ will be assumed. Then
$$j_v(y)\otimes j_u(x)=j_{u\times v}(x\otimes y)\qquad (u\times
v)^c=v^c\times u^c$$
for each $u,v\in\R{G}$.

For each $f\in\Mor(u,v)$ we shall denote by $f^c\in
\Mor(\bar{u},\bar{v})$ a morphism specified by
\begin{equation}\label{cmor}
f^cj_u=j_vf.
\end{equation}

To every representation $u\in\R{G}$ we can associate two elementary
intertwiners $\I_u\in\Mor(\1,u\bar{u})$ and $\btwn^u\in
\Mor(\bar{u}u,\1)$, given by the identity operator and the contraction
map of $H_u$ respectively. Because of the symmetry between $u$ and
$\bar{u}$ we have also intertwiners $\I^u\in\Mor(\1,\bar{u}u)$ and
$\btwn_u\in\Mor(u\bar{u},\1)$. In explicit form,
\begin{equation*}
\begin{gathered}
\btwn^u(f\otimes x)=f(x)\\
\I_u(1)=\sum_i e_i\otimes e_i^*
\end{gathered}\qquad
\begin{gathered}
\btwn_u(x\otimes f)=fC_u(x)\\
\I^u(1)=\Sum_{ij}[C_u^{-1}]_{ji} e_i^*\otimes e_j,
\end{gathered}
\end{equation*}
where $\{e_i\}$ is an arbitrary orthonormal basis in $H_u$ and
$\{e_i^*\}$ is the corresponding biorthogonal basis.

We shall use the symbol $\otimes_M$ for the tensor product over
$\cal{V}$, in order to stress the geometrical background for our
considerations.

\begin{pro} The following natural bimodule isomorphism holds
\begin{equation}
\bim{u\times v}\leftrightarrow \bim{u}\otimes_M\!\bim{v}.\label{fun-prod}
\end{equation}
The above isomorphism is induced by the product map of $\cal{B}$.
\end{pro}
\begin{pf}
Let us fix a representation $u$. As follows from \cite{D2}-Appendix A,
there exist intertwiners $\mu_k,\nu_k\in\bim{\bar{u}},\bim{u}$ such
that
$$\sum_k \mu_k(f)\nu_k(x)=f(x)1,$$
for each $f\in H_{\bar{u}}$ and $x\in H_u$.

Let $\vartheta_{vu}\colon\bim{v\times u}\rightarrow
\bim{v}\otimes_M\!\bim{u}$ be a map given by
$$ \vartheta_{vu}(\xi)=\sum_{k}\xi_k \otimes\nu_k,$$
where $\xi_k(f)=\Sum_{j}\xi(f\otimes
e_j)\mu_k(e_j^*)$.

We prove that the inverse of $\vartheta_{vu}$ is given by
$\vartheta_{vu}^{-1}(\varphi\otimes\psi)(y\otimes x)=
\varphi(y)\psi(x)$. A direct calculation gives
$$\vartheta_{vu}^{-1}\vartheta_{vu}(\xi)(y\otimes x)=
\sum_{kj}\xi(y\otimes
e_j)\mu_k(e_j^*)\nu_k(x)=\xi(y\otimes x), $$
and similarly
$$
\vartheta_{vu}\vartheta^{-1}_{vu}(\varphi\otimes \psi)
=\sum_{kj}
\varphi\bigl\{\psi(e_j)\mu_k(e_j^*)\bigr\}\otimes\nu_k=\Sum_{kj}
\varphi\otimes \psi(e_j)\mu_k(e_j^*)\nu_k=\varphi\otimes \psi,
$$
and hence $\vartheta_{vu}$ is bijective.
\end{pf}

In terms of the above identification, the elementary intertwiners $\I_u$
and $\btwn^u$ induce bimodule homomorphisms
$$\pre{\btwn^u}\colon\cal{V}\rightarrow\bim{\bar{u}}\otimes_M\!\bim{u}
\qquad
\langle\rangle_u^+=\pre{[\I_u]}\colon\bim{u}\otimes_M\!\bim{\bar{u}}
\rightarrow\cal{V}$$
respectively. Explicitly, it follows that
\begin{equation}
\pre{\btwn^u}(1)=\sum_k\mu_k\otimes\nu_k\qquad
\langle\varphi\otimes\psi\rangle_u^+=\sum_i\varphi(e_i)\psi(e_i^*).
\end{equation}

Interchanging the roles of $u$ and $\bar{u}$ we obtain a
pairing $\langle\rangle_u^-
\colon\bim{\bar{u}}\otimes_M\!\bim{u}\rightarrow\cal{V}$ and the
injection $\btwn_u^\sstar\colon\cal{V}\rightarrow\bim{u}\otimes_M\!
\bim{\bar{u}}$.

Maps $j_u$ naturally induce antilinear antiisomorphisms
$\amap{u}\colon\bim{u}\rightarrow\bim{\bar{u}}$,
via the following diagram
\begin{equation}\label{amap}
\begin{CD}
H_{u} @>{\mbox{$\phantom{\amap{u}}\varphi$}}>> \cal{B}\\
@V{\mbox{$j_u$}}VV @VV{\mbox{$*$}}V\\
H_{u}^* @>>{\mbox{$\amap{u}\varphi$}}> \cal{B}
\end{CD}
\end{equation}
With the help of $\amap{u}$,
the bimodule $\bim{\bar{u}}$ is naturally identificable with
the {\it conjugate} of $\bim{u}$.

\begin{lem}
The following identities hold
\begin{gather*}
\amap{\bar{u}}\amap{u}=\pre{[C^{-1}_u]}\\
\pre{f^c}=\amap{u}\pre{f}\amap{v}^{-1}\\
\amap{v}(\psi)\otimes\amap{u}(\varphi)=\amap{u\times v}(\varphi\otimes
\psi),
\end{gather*}
where $f\in\Mor(u,v)$ and $\varphi,\psi\in\bim{u,v}$. \qed
\end{lem}

Let us observe that the first two equalities allow us to extend the
introduced representation of morphisms, to all {\it antilinear
morphisms}. By definition, an antilinear map
$f\colon H_u\rightarrow H_v$ is a {\it morphism}, iff the diagram
\begin{equation}
\begin{CD}
H_u @>{\mbox{$u$}}>> H_u\otimes\cal{A}\\
@V{\mbox{$f$}}VV @VV{\mbox{$f\otimes *$}}V\\
H_v @>>{\mbox{$v$}}> H_v\otimes\cal{A}
\end{CD}
\end{equation}
is commutative.  Let us denote by $\M(u,v)$ the corresponding
spaces of antilinear morphisms. From this moment, we shall include
all such maps in the system of morphisms of $\R{G}$.
Consequently, let us assume that
$$\MoR(u,v)=\Mor(u,v)\oplus\M(u,v) $$
represent arrows of $\R{G}$. These spaces will be endowed with their
natural $\Bbb{Z}_2$-grading.

Let $\Phi_M$ be the category of finite projective $\cal{V}$-bimodules.
For each $\cal{J},\cal{K}\in\Phi_M$, let
$\H(\cal{J},\cal{K})$ be the space of the corresponding
antilinear bimodule antihomomorphisms, and let us consider the spaces
$$\HoM(\cal{J},\cal{K})=\Hom(\cal{J},\cal{K})\oplus\H(\cal{J},\cal{K}),$$
endowed with the natural $\Bbb{Z}_2$-grading. We shall assume that
the arrows of $\Phi_M$ are given by the elements of
$\HoM(\cal{J},\cal{K})$.

If $f\in\M(u,v)$ then the formula
\begin{equation}
\pre{f}=\pre{[j_vf]}\amap{v}
\end{equation}
defines a grade-preserving
linear extension $\map{uv}\colon\MoR(u,v)
\rightarrow\HoM(\bim{u},\bim{v})$
of the previously introduced representation. The formula
\eqref{contr} remains valid for arbitrary $f,g\in\MoR(u,v)$.

Let $\cal{T}\subseteq\R{G}$ be a complete set of mutually inequivalent
irreducible unitary representations of $G$.
The algebra $\cal{B}$ can be decomposed into a direct sum of multiple
irreducible $\cal{V}$-bimodules
$$\cal{B}=\sideset{}{^\oplus}\sum_{\alpha\in\cal{T}}\cal{B}^\alpha,$$
relative to the right action $F$. Each $\cal{B}^\alpha$ can be further
decomposed as
$$\cal{B}^\alpha=\Mor(\alpha,F)\otimes H_\alpha\qquad
\varphi(x)\leftrightarrow\varphi\otimes x.$$
In this sense
$\bim{\alpha}=\Mor(\alpha,F)$ are elementary building
blocks of the bundle $P$.

Motivated by the above considerations, we are going to formulate
a notion of a representation of a concrete monoidal $W$-category.

Let $R\subseteq\R{G}$ be a generating monoidal subcategory, closed (up
to the equivalence) under the conjugation functor and containing all the
conjugation maps.

\smallskip
\bla{rep1} Let us assume that we have a pair
$$\varrho=\Bigl(\Bigl\{\bim{u};u\in R\Bigr\}
\Bigr\{\map{uv};u,v\in R\Bigr\}\Bigr)$$
where $\bim{u}$ are finite projective $\cal{V}$-bimodules,
$\map{uv}\colon\MoR(u,v)\rightarrow\HoM(\bim{v},\bim{u})$ are
grade-preserving linear maps,
realizing a contravariant functor $\varrho\colon R\rightarrow\Phi_M$. We
shall also use the simplified notation $\map{uv}(f)=\pre{f}$.

\smallskip
Let us consider the category $R\times R$  defined
in the following way. The objects are the corresponding ordered pairs,
while the morphisms between $(u,p)$ and $(v,q)$ are pairs
$(\varphi,\psi)$ satisfying
$\varphi\in\Mor(u,v)$ and $\psi\in\Mor(p,q)$ for the degree $0$, and
satisfying $\varphi\in\M(u,q)$ and $\psi\in\M(p,v)$ for the degree $1$.
In a similar way, let us introduce the morphisms in
$\Phi_M\times\Phi_M$. Let us denote by $[\varphi,\psi]$ the twisted
tensor product of the corresponding antilinear transformations.
The operations $\times$ and $\otimes_M$ are
undersandable as grade-preserving covariant functors
$$\times\colon R\times R\rightarrow R\qquad {\otimes_M}\colon
\Phi_M\times\Phi_M\rightarrow\Phi_M,$$
transforming $(\varphi,\psi)$ into $\varphi\otimes\psi$ or
$[\varphi,\psi]$, depending on the parity.

\smallskip
\bla{rep2}
Furthermore, let us assume that there exists a natural transformation
$\vartheta$ between functors
$\varrho\circ\{{\times}\}$ and $\otimes_M\!\circ(\varrho\times\varrho)$,
given by the system of bimodule isomorphisms
$\vartheta_{uv}\colon\bim{u\times
v}\rightarrow\bim{u}\otimes_M\!\bim{v}$.

\smallskip
\bla{rep3}
Finally, let us assume that $\vartheta$
is associative, in the sense that
\begin{equation}
(\vartheta_{uv}\otimes\id)\vartheta_{u\times
v,w}=(\id\otimes\vartheta_{vw})\vartheta_{u,v\times w},
\end{equation}
for each $u,v,w\in R$.

\begin{defn} Every pair $\tau=(\varrho,\vartheta)$ satisfying the
above listed
properties is called a {\it representation} of $R$ in $\Phi_M$.
\end{defn}

The analysis of this section may be now summarized as follows
\begin{pro}
Let $P=(\cal{B},i,F)$ be an arbitrary quantum principal
$G$-bundle over $M$. Then the corresponding intertwiner bimodules
$\bim{u}$, together with the associated system of maps $\map{uv}$
and the product identifications form a
representation $\tau_P$ of $\R{G}$ in $\Phi_M$.\qed
\end{pro}

\section{Bundle Reconstruction}
The aim of this section is to prove that the construction of the
previous section works in both directions, so that there exists a
natural bijection
$$\Bigl\{\,\mbox{Quantum $G$-bundles $P$ over $M$}\,\Bigr\}
\,\leftrightarrow\,\Bigl\{\,\mbox{Representations $\tau$
of $\R{G}$ in $\Phi_M$}\,\Bigr\}.$$
Let $R$ be a generating subcategory of $\R{G}$, closed (up to the
equivalence) under taking conjugate objects.
\begin{pro} Every representation of $R$ in $\Phi_M$ can be
extended to a representation of $\R{G}$. The extension is unique, up to
a natural transformation.
\end{pro}
\begin{pf} Let $\tau$ be a representation of $R$ in $\Phi_M$. Let
us first assume that the extension exists, and
prove its uniqueness. Each object $u\in\R{G}$ can be realized as a direct
summand in some $r=\Sum_k^\oplus r_k$, where $r_k\in R$, with the help
of embedding and projection morphisms, $\iota_u\in\Mor(u,r)$ and
$\pi_u\in\Mor(r,u)$. Then the object $\bim{u}$ is realizable as a
submodule in $\bim{r}=\Sum^\oplus_k\bim{r_k}$, with the help of the
embedding $\pre{[\pi_u]}$ and the projection $\pre{[\iota_u]}$.

Let us consider a representation $v$, realized in a similar
way in $s=\Sum_k^\oplus s_k$. For each $f\in\Mor(u,v)$ the bimodule
homomorphism $\pre{f}$ is uniquely fixed, because we can
write $f\leftrightarrow\Sum_{kl} i_l\rho_{kl}p_k$,
where $i_l\in\Mor(s_l,s)$ and
$p_k\in\Mor(r,r_k)$ are canonical embedings and projections, and
$\rho_{kl}\in\Mor(r_k,s_l)$. This means that all
maps $\pre{f}$ are expressible in terms
of $R$. The same holds for antilinear morphisms.
Hence, the extension is unique.

Conversely, starting from the above observations it is possible to
{\it define} the extension of $\tau$. Objects and morphisms of
$\R{G}$ are expressible purely in terms
of $R$, and such expressions consistently define the extension of $\tau$.
In particular, bimodules $\bim{u}$
can be invariantly described by ``gluing'' the corresponding images
in various possible $\bim{r}$. We omit technical details.
\end{pf}

Let $\tau$ be an arbitrary representation of $\R{G}$ in
$\Phi_M$. We are going to re-construct the bundle $P$ which satisfies
$\tau=\tau_P$.

As first, let us define a $\cal{V}$-bimodule $\cal{B}$ as a direct sum
\begin{equation}\label{def-B}
\cal{B}=\sideset{}{^\oplus}\sum_{\alpha\in\cal{T}}\bim{\alpha}\otimes
H_\alpha.
\end{equation}
The group $G$ naturally acts on $\cal{B}$ on the right, via the sum $F$
of actions $\id\otimes \alpha$. Let $i\colon\cal{V}\rightarrow\cal{B}$
be the canonical inclusion map. By definition, $b\in\im(i)$ iff
$F(b)=b\otimes 1$ for each $b\in\cal{B}$.

Let us further observe that the following natural isomorphism holds
\begin{equation}\label{Mor-E}
\bim{u}\leftrightarrow\Mor(u,F),
\end{equation}
for each $u\in\R{G}$. Indeed, if $u=\alpha\in\cal{T}$ then the above
correspondence is simply $\varphi\colon x\mapsto\varphi\otimes x$.

Let us consider an arbitrary representation $u\in\R{G}$. The following
intrinsic bimodule decomposition holds
\begin{equation}\label{dec1}
\bim{u}=\sideset{}{^\oplus}\sum_{\alpha\in\cal{T}}\Mor(u,\alpha)\otimes
\bim{\alpha},
\end{equation}
where actually the sum is finite. On the other hand
\begin{equation}\label{dec2}
\Mor(u,F)=\sideset{}{^\oplus}\sum_{\alpha\in\cal{T}}\Mor(u,\alpha)\otimes
\Mor(\alpha,F),
\end{equation}
where the identification is induced by the composition of morphisms.
Combining the above decompositions, we conclude
that \eqref{Mor-E} holds in the full generality.
In terms of this identification,
the following correspondence holds
\begin{equation}
\pre{f}(\varphi)\leftrightarrow\varphi f,
\end{equation}
for each $u,v\in\R{G}$, $\varphi\in\bim{v}$ and $f\in\Mor(u,v)$.
In what follows we shall assume all the above identifications.

For $\varphi\in\bim{u}$ and $\psi\in\bim{v}$ we shall denote by
$\varphi\psi\in\bim{u\times v}$ the element given by
\begin{equation}
\varphi\psi=\vartheta_{uv}^{-1}(\varphi\otimes\psi).
\end{equation}
The associativity of the natural
transformation $\vartheta$ implies that the above product is
associative, too.

\begin{pro} The formula
\begin{equation}\label{prod-B}
\varphi(x)\psi(y)=(\varphi\psi)(x\otimes y)
\end{equation}
where $\varphi\in\bim{u}$ and $\psi\in\bim{v}$, consistently defines a
structure of a unital associative algebra in $\cal{B}$ so that maps
$\{\,i,F\,\}$ are unital homomorphisms.
\end{pro}
\begin{pf}
The above formula, restricted to representations from $\cal{T}$, defines
a bilinear product in $\cal{B}$, so
that the embedding $i$ is multiplicative. The fact that this product
satisfies the same formula for arbitrary representations follows from
the naturality of $\vartheta$. Further, $i(1)=1\otimes 1$ is
the unit element. In particular $F$ is unital, too. The map $F$
preserves the introduced product because
\begin{multline*}
F\bigl(\varphi(e_i)\psi(e_j)\bigr)=
F\bigl((\varphi\psi)(e_i\otimes e_j)\bigr)=
\sum_{kl}(\varphi\psi)(e_k\otimes e_l)\otimes u_{ki}
v_{lj}\\=
\sum_{kl}\varphi(e_k)\psi(e_l)\otimes u_{ki}
v_{lj}
=\sum_{kl}\bigl(\varphi(e_k)\otimes u_{ki}\bigr)\bigl(
\psi(e_l)\otimes v_{lj}\bigr)=F\varphi(e_i)F\psi (e_j).
\end{multline*}
As a consequence of the associativity of the product of intertwiners we
obtain
\begin{multline*}
\varphi(x)\bigl(\psi(y)\eta(z)\bigr)=\varphi
(x)\bigl((\psi\eta)(y\otimes z)\bigr)=
\bigl(\varphi(\psi\eta)\bigr)(x\otimes y\otimes
z)\\=\bigl((\varphi\psi)\eta\bigr)(x\otimes
y\otimes z)=(\varphi\psi)(x\otimes
y)\eta(z)=\bigl(\varphi(x)\psi(y)\bigr)
\eta(z),
\end{multline*}
which completes the proof.
\end{pf}

Now we shall define the *-structure on the algebra $\cal{B}$. There
exists the unique antilinear map $*\colon\cal{B}\rightarrow\cal{B}$
satisfying
$$[\pre{f}(\psi)(x)]^*=\psi[f(x)]$$
for each $f\in\M(v,u)$ and $\psi\in\bim{u}$. The consistency of
this definition directly follows from the functoriality of $\tau$.

\begin{pro}
The introduced map is a *-structure on the algebra $\cal{B}$.
\end{pro}
\begin{pf}
Let $\amap{u}\colon\bim{u}\rightarrow\bim{\bar{u}}$
be antilinear bimodule
homomorphisms given by $$\amap{u}=\pre{[j_u^{-1}]}.$$
We have then
$$
\psi(x)^{**}=\bigl[(\amap{u}\psi)(j_ux)\bigr]^*=
(\amap{\bar{u}}\amap{u}\psi)(j_{\bar{u}}j_ux)
=\bigl\{\pre{[C_u^{-1}]}\psi\bigr\}(C_ux)=\psi(x),
$$
which shows that $*$ is involutive. Furthermore, if $\varphi\in\bim{u}$
and $\psi\in\bim{v}$ then
\begin{multline*}
\bigl(\varphi(x)\psi(y)\bigr)^*=\bigl[
\amap{u\times v}(\varphi\psi)\bigr](j_v(y)\otimes
j_u(x))\\=
\bigl(\amap{v}\psi\bigr)\bigl(j_v(y)\bigr)\bigl(
\amap{u}\varphi\bigr)\bigl(j_u(x)\bigr)
=\psi(y)^*\varphi(x)^*,
\end{multline*}
in accordance with the naturality of the transformation $\vartheta$.
Hence, $*$ is antimultiplicative.
\end{pf}

Summarizing the analysis of this section, we conclude that
\begin{thm} The triplet $P=(\cal{B},i,F)$ is a quantum principal
$G$-bundle over $M$. Moreover, $\tau=\tau_P$.
\end{thm}
\begin{pf}
It remains to prove that $G$ is acting freely on $P$. For a given
$u\in\R{G}$ let us consider elements $\mu_k\in\bim{\bar{u}}$ and
$\nu_k\in\bim{u}$ such that
$$
\pre{\btwn^u}(1)\leftrightarrow\Sum_k\mu_k\otimes \nu_k.
$$
We have then
$$\sum_k\mu_kj_u(x)\nu_k(y)=(x,y)1,$$
for each $x,y\in H_u$. This implies that $F$ is free.
\end{pf}
\section{Examples}
\subsection{Bundles Over $1$-point Sets}
If $M$ is a $1$-point set then $\cal{V}=\Bbb{C}$ and objects of
$\Phi_M$ are complex finite-dimensional vector spaces.
Morphisms are linear maps, and binary operations are standard direct
sums and tensor products of spaces. The problem of classifying
$G$-bundles over such a trivial space $M$ thus reduces to a
linear-algebraic game. Interestingly, all non-trivial examples of
bundles of this type will be completely quantum--without points at all.

\subsection{Quantum Line Bundles}
Let us assume that $G=U(1)$. Then $\cal{A}$ is generated by a single
unitarity $u$ satisfying $\phi(u)=u\otimes u$, identified with the
fundamental representation. Irreducible representations are labeled by
integers, $u^n\leftrightarrow n$. Category $\R{G}$ is generated by
objects $\bar{u}=u^{-1}$ and $u$. The only relations in $\R{G}$ are
identifications
$$\bar{u}\times u\leftrightarrow\1\qquad u\times
\bar{u}\leftrightarrow\1,$$
together with the compatibility conditions
\begin{gather*}
(\bar{u}\times u)\times \bar{u}\leftrightarrow \1\times
\bar{u}=\bar{u}\times\1\leftrightarrow \bar{u}\times (u\times
\bar{u})\\
(u\times \bar{u})\times u\leftrightarrow \1\times
u=u\times\1\leftrightarrow u\times (\bar{u}\times u).
\end{gather*}

Let $\tau$ be a representation of $\R{G}$ in $\Phi_M$. It follows that
$\tau$ is completely specified by a
$\cal{V}$-bimodule $\cal{E}$ and its conjugate $\bar{\cal{E}}$,
together with hermitian bimodule isomorphisms
$\bar{\mu}\colon\cal{E}\otimes_M\!\bar{\cal{E}}\rightarrow\cal{V}$
and $\mu\colon\bar{\cal{E}}\otimes_M\!\cal{E}
\rightarrow\cal{V}$, which are mutually compatible such that
\begin{equation*}
\mu\otimes\id=\id\otimes\bar{\mu}\qquad\bar{\mu}\otimes\id=\id\otimes
\mu.
\end{equation*}

As a concrete illustration, let us assume that $M$ is a classical
compact smooth manifold, and $\cal{V}=S(M)$. Then, at the level of left
modules, the elements of $\cal{E}$ can be naturally
identified with smooth sections of a line bundle $L$ over $M$. In
terms of this identification the right $\cal{V}$-module structure is
determined by a *-automorphism
$\varepsilon\colon S(M)\rightarrow S(M)$, so that $\varphi
f=\varepsilon(f)\varphi$, for each $\varphi\in\cal{E}$ and $f\in S(M)$.
The presence of $\varepsilon$ does not influence the bimodule structure
of $\cal{E}\otimes_M\!\bar{\cal{E}}$ and the vector space
structure of
$\bar{\cal{E}}\otimes_M\!\cal{E}$. Hence it is possible to introduce
classical natural maps
$\mu_0\colon\bar{\cal{E}}\otimes_M\!\cal{E}\rightarrow\cal{V}$
and $\bar{\mu}_0\colon\cal{E}\otimes_M\!\bar{\cal{E}}
\rightarrow\cal{V}$. We have
$$
\mu_0(f\psi g)=\varepsilon(f)\mu_0(\psi)
\varepsilon(g)\qquad\bar{\mu}_0(f\varphi g)=f\bar{\mu}_0(\varphi)g
$$
for each $\psi\in\bar{\cal{E}}\otimes_M\!\cal{E}$, $\varphi\in
\cal{E}\otimes_M\!\bar{\cal{E}}$ and $f,g\in\cal{V}$. This implies
that there exist invertible smooth functions $U,V\colon
M\rightarrow\Bbb{C}$ satisfying
$$\mu(\psi)=U\varepsilon^{-1}\mu_0(\psi)\qquad\bar{\mu}(\varphi)=
V\bar{\mu}_0(\varphi).$$
Finally, compatibility conditions between $\mu$ and $\bar{\mu}$ imply
$$V=\varepsilon(U)\qquad U=U^*.$$

Therefore, quantum line bundles over $M$ are classified by triplets
$(L,\varepsilon,U)\leftrightarrow P$ of the above described type.

\subsection{Universal Unitary Groups}

Let us assume that $G$ is a {\it universal} \cite{D1}-Appendix~A
compact matrix quantum
group, with the fundamental representation $u$. Let
$\cal{C}$ be a concrete monoidal $W$-category \cite{W2} generated by
elements $\{u,\bar{u}\}$. The objects of $\cal{C}$ are
just the words over $\{u,\bar{u}\}$, including the unit object.
The morphisms between objects of $\cal{C}$ are generated by elementary
morphisms $\{\btwn^u, \btwn_u,\I_u,\I^u\}$, and the conjugation maps
$\{j_u, j_u^{-1}\}$. The only relations between standard morphisms in
$\cal{C}$ are given by
\begin{gather*}
(\btwn^u\otimes\id)(\id\otimes
\I_u)=\id=(\id\otimes\btwn_u)(\I^u\otimes\id)\\
(\btwn_u\otimes\id)(\id\otimes
\I^u)=\id=(\id\otimes\btwn^u)(\I_u\otimes\id)\\
\btwn_u \I_u=\btwn^u \I^u=n_u,
\end{gather*}
where $n_u=\tr(C_u)=\tr(C_u^{-1})$. We have also the hermicity
conditions
\begin{gather*}
\btwn_u[j_u,j_u^{-1}]=\btwn_u\qquad
\btwn^u[j^{-1}_u,j_u]=\btwn^u\\
[j_u,j_u^{-1}]\I_u=\I_u\qquad[j^{-1}_u,j_u]\I^u=\I^u.
\end{gather*}

Therefore, the representations
$\tau$ of $\cal{C}$ are labeled by conjugate $\cal{V}$-bimodules
$\bim{u},\bim{\bar{u}}$, together with hermitian bimodule injections
$\btwn_u^\sstar$ and $\pre{\btwn^u}$ and contractions
$\langle\rangle_u^\pm$, reflecting the above relations in $\cal{C}$. In
particular, the nature of the structure group influences
the whole classification only via the single positive number $n_u$.

\subsection{Some Variations}

The presented formalism can be applied to the study of differential
structures on quantum principal bundles. Let us assume that the full
calculus on $P$ is described by a graded-differential algebra
$\Omega(P)$, and let $\hor(P)\subseteq\Omega(P)$ be the *-subalgebra
representing horizontal forms \cite{D2}. Let
$F^\wedge\colon\hor(P)\rightarrow\hor(P)\otimes\cal{A}$ be the right
action map. Let $\Omega(M)$ be a graded-differential *-algebra
representing the calculus on $M$. It consists of $F^\wedge$-invariant
elements.

Let $\tau^\star$ be the corresponding representation of $\R{G}$ in
the category $\Phi^\sstar_M$ of finite projective $\Omega(M)$-bimodules,
with intertwiner bimodules $\bimh{u}$. These bimodules are
naturally graded, and $\bimh{u}^0=\bim{u}$ for each $u\in\bim{u}$.
Furthermore, the following natural decompositions hold
\begin{equation}\label{WM-M}
\bim{u}\otimes_M\!\Omega(M)\leftrightarrow\bimh{u}\leftrightarrow
\Omega(M)\otimes_M\!\bim{u}.
\end{equation}
These decompositions are induced by the left/right product maps in
$\bimh{u}$.

The composition of the two identifications defines a grade-preserving
flip-over operator $$\sigma_u\colon\bim{u}\otimes_M\!\Omega(M)
\rightarrow\Omega(M)\otimes_M\!\bim{u},$$
acting as the identity on $\bim{u}$. This operator contains the whole
information concerning the $\Omega(M)$-bimodule structure.

The following compatibility condition expresses
the consistency of the $\Omega(M)$-bimodule structure. Let
$m_M$ be the product in $\Omega(M)$. It follows that
\begin{equation}
\sigma_u(\id\otimes m_M)=(m_M\otimes\id)
(\id\otimes\sigma_u)(\sigma_u\otimes\id).
\end{equation}
A similar compatibility condition holds between operators
$\sigma_u$ and the natural transformation $\vartheta$. We have
\begin{equation}
(\id\otimes\vartheta_{uv})
\sigma_{u\times
v}=(\sigma_u\otimes\id)(\id\otimes\sigma_v)(\vartheta_{uv}\otimes\id)
\end{equation}
for each $u,v\in\R{G}$.

Finally, intertwiner homomorphisms
$\pre{f}$ and conjugation maps $\amap{u}$ satisfy the following
diagrams
\begin{equation*}
\begin{CD} \bim{v}\otimes_M\!\Omega(M) @>{\mbox{$\sigma_v$}}>>
\Omega(M)\otimes_M\!\bim{v}\\
@V{\mbox{$\pre{f}\otimes\id$}}VV @VV{\mbox{$\id\otimes\pre{f}$}}V\\
\bim{u}\otimes_M\!\Omega(M) @>>{\mbox{$\sigma_u$}}>
\Omega(M)\otimes_M\!\bim{u}
\end{CD}\qquad
\begin{CD} \bim{u}\otimes_M\!\Omega(M) @>{\mbox{$\sigma_u$}}>>
\Omega(M)\otimes_M\!\bim{u}\\
@V{\mbox{$[\amap{u},*]$}}VV  @VV{\mbox{$[*,\amap{u}]$}}V\\
\Omega(M)\otimes_M\!\bim{\bar{u}} @<<{\mbox{$\sigma_{\bar{u}}$}}<
\bim{\bar{u}}\otimes_M\!\Omega(M)
\end{CD}
\end{equation*}

The algebra of horizontal forms can be decomposed in the following way
\begin{equation}
\begin{gathered}
\hor(P)=\sideset{}{^\oplus}\sum_{\alpha\in\cal{T}}\cal{H}^\alpha(P)\qquad
\cal{H}^\alpha(P)=\bimh{\alpha}\otimes H_\alpha\\
\hor(P)\leftrightarrow\Omega(M)\otimes_M\!\cal{B}
\leftrightarrow\cal{B}\otimes_M\!\Omega(M).
\end{gathered}
\end{equation}

Let us assume that $P$ admits regular connections \cite{D2}. For a given
regular connection $\omega$, let $D\colon\hor(P)\rightarrow\hor(P)$ be
its covariant derivative. This map is a right-covariant hermitian
first-order antiderivation.  In particular, it naturally
induces (via compositions) first-order
maps $D_u\colon\bimh{u}\rightarrow\bimh{u}$. The following properties hold
\begin{gather*}
D_{u\times v}(\varphi\otimes\psi)=D_u(\varphi)\psi
+(-1)^{\partial\varphi}\varphi D_v(\psi)\qquad D_{\1}=d_M\\
D_{\bar{u}}\amap{u}=\amap{u}D_{u}\qquad D_u\map{uv}(f)=\map{uv}(f)D_v
\end{gather*}
for each $u,v\in\R{G}$ and $f\in\Mor(u,v)$. Here $d_M$ is the
differential on $\Omega(M)$.

Conversely, starting from a representation $\tau^\star$ (with graded
bimodules, and grade-preserving maps between them) we can reconstruct
$\hor(P)$. If we start from $\hor(P)$, then the construction
of the whole differential calculus on $P$ can be completed applying
methods presented in \cite{D3}. Covariant derivatives of
regular connections are in a natural correspondence
with systems
$\Bigl\{D_u;u\in\R{G}\Bigr\}$ of linear maps
$D_u\colon\bimh{u}\rightarrow\bimh{u}$ satisfying the above conditions.

\end{document}